\documentclass[11pt,leqno]{amsart}

\usepackage{amssymb,multirow}

\newcommand{\leb}{\ensuremath{l^1}}
\newcommand{\oT}{\ensuremath{^{\{1\}}}}
\newcommand{\tT}{\ensuremath{^{\{2\}}}}
\newcommand{\jT}{\ensuremath{^{\{j\}}}}
\newcommand{\re}{\ensuremath{\textrm{Re}}}

\newtheorem{theorem}{Theorem}[section]
\newtheorem{lemma}[theorem]{Lemma}
\newtheorem{proposition}[theorem]{Proposition}

\theoremstyle{remark}
\newtheorem{remark}[theorem]{Remark}

\numberwithin{equation}{section}

\title[Large data local well-posedness]{Large data local well-posedness for a class of KdV-type equations}
\author[B. Harrop-Griffiths]{Benjamin Harrop-Griffiths}
\address{Department of Mathematics, University of California, Berkeley, CA 94720}
\email{benhg@math.berkeley.edu}

\subjclass[2010]{Primary 35Q53, 35G25}

\begin{document}

\begin{abstract}
 In this article we consider the Cauchy problem with large initial data for an equation of the form \[(\partial_t+\partial_x^3)u=F(u,u_x,u_{xx})\]where \(F\) is a polynomial with no constant or linear terms. Local well-posedness was established in weighted Sobolev spaces by Kenig-Ponce-Vega. In this paper we prove local well-posedness in a translation invariant subspace of \(H^s\) by adapting the result of Marzuola-Metcalfe-Tataru on quasilinear Schr\"odinger equations.
\end{abstract}

\maketitle

\thispagestyle{empty}

\section{Introduction}

In this article we consider local well-posedness for an equation of the form

\begin{equation}\label{eq:nonlinear}
\left\{\begin{array}{l}(\partial_t+\partial_x^3)u=F(u,u_x,u_{xx}),\qquad u\colon\mathbb{R}\times\mathbb{R}\rightarrow\mathbb{R}\textrm{ or }\mathbb{C}\vspace{3pt}\\u(0)=u_0\end{array}\right.
\end{equation}
where we assume \(F\) is a constant coefficient polynomial of degree \(m\geq 2\) with no linear or constant terms.

It is natural to consider well-posedness in \(H^s(\mathbb{R})\). However, due to the infinite speed of propagation, even a linear equation of the form
\[
(\partial_t+\partial_x^3+a(x)\partial_x^2)u=0
\]
where \(a\) is smooth with bounded derivatives requires a Mizohata-type necessary condition for \(L^2\) well-posedness \cite{A,Mz,Tam}
\begin{equation}\label{eq:mizohata}
\sup_{x_1\leq x_2}\re\int_{x_1}^{x_2} a(x)\,dx<\infty
\end{equation}
So at the very least, if \(F\) contains a term of the form \(uu_{xx}\), then we expect any solution \(u\) to \eqref{eq:nonlinear} to require some additional integrability. Indeed, an ill-posedness result in \(H^s\) was proved by Pilod \cite{P}.

One way to address this difficulty is to consider weighted spaces. Kenig-Ponce-Vega proved local well-posedness for small data in \cite{KPV1} and arbitrary data in \cite{KPV2} using the weighted space \(H^s\cap L^2(|x|^k\,dx)\) for sufficiently large \(k\in\mathbb{Z}^+\) and \(s>0\). Their result was extended to systems by Kenig-Staffilani \cite{KS}. Replacing weighted \(L^2\) spaces with weighted Besov spaces, Pilod \cite{P} proved local well-posedness for small data at low regularities in the space \(H^s(\mathbb{R})\cap\mathcal{B}^{s-2,2}_2(\mathbb{R},x^2dx)\) where \(s>\tfrac{9}{4}\) for certain quadratic nonlinearities. Quasilinear versions of this problem for which \eqref{eq:nonlinear} is a special case have also been studied by several authors (see \cite{A,C} and references therein).

As the equation \eqref{eq:nonlinear} is translation invariant, it is natural to look for a solution in a translation invariant space. By replacing weighted spaces with a spatial summability condition, Marzuola-Metcalfe-Tataru \cite{MMT} proved a small data result for quasilinear Schr\"odinger equations for initial data in a translation invariant space \(l^1H^s\subset H^s\). Their result relies on a local energy decay estimate using spaces similar to those suggested by Kenig-Ponce-Vega \cite{KPV3}.

In this paper we adapt the result of Marzuola-Metcalfe-Tataru to the problem \eqref{eq:nonlinear} where waves at frequency \(2^j\) travel at speed \(2^{2j}\). By a slight abuse of notation we also call the adapted initial data space \(\leb H^s\). As the need for additional integrability is solely due to the bilinear interactions, as in \cite{KPV1,KPV2,KPV3,KS,MMT}, we expect to be able to remove the spatial summability condition for the case that \(F\) contains no quadratic terms.

We take a standard Littlewood-Paley decomposition
\[
1=\sum\limits_{j=0}^\infty S_j
\]
constructed by taking smooth \(\varphi_0\colon\mathbb{R}\rightarrow[0,1]\) such that
\[
\varphi_0(\xi)=\left\{\begin{array}{ll}1&\qquad\textrm{for }\xi\in[-1,1]\vspace{3pt}\\0&\qquad\textrm{for }|\xi|\geq2\end{array}\right.
\]
Then define
\[
\varphi_j(\xi)=\varphi_0(2^{-j}\xi)-\varphi_0(2^{-j+1}\xi)
\]
and
\[
f_j=S_jf=\mathcal{F}^{-1}(\varphi_j\hat f)
\]
where \(\mathcal{F}u=\hat u\) is the spatial Fourier transform.

For each \(j\geq0\) we take a partition \(\mathcal{Q}_{2j}\) of \(\mathbb{R}\) into intervals of length \(2^{2j}\) and an associated smooth partition of unity
\[
1=\sum\limits_{Q\in\mathcal{Q}_{2j}}\chi_Q
\]
where we assume \(\chi_Q\sim 1\) on \(Q\) and \(\mathrm{supp}\,\chi_Q\subset B\left(Q,\tfrac{1}{2}\right)\) then define
\[
\|u\|_{l^1_{2j}L^2}=\sum\limits_{Q\in\mathcal{Q}_{2j}}\|\chi_Q u\|_{L^2}
\]
We define the initial data space \(\leb H^s\) with norm
\[
\|u\|_{\leb H^s}^2=\sum\limits_{j\geq0}2^{2sj}\|S_ju\|^2_{l^1_{2j}L^2}
\]
We note that for \(s>1\) we have \(\leb H^s\subset L^1\). The main result we prove is the following.

\vspace{15pt}\begin{theorem}\label{mainthrm}
For \(s>\tfrac{9}{2}\), there exists \(C>0\) such that the equation \eqref{eq:nonlinear} is locally well-posed in \(\leb H^s\) on the time interval \([0,T]\) where \(T=e^{-C\|u_0\|_{\leb H^s}}\).
\end{theorem}\vspace{15pt}

We take the definition of ``well-posedness'' to be the existence and uniqueness of a solution \(u\in C([0,T],\leb H^s)\) and Lipschitz continuity of the solution map\[\leb H^s\ni u_0\mapsto u\in C([0,T],\leb H^s)\]We note that as the problem considered in \cite{MMT} is quasilinear, continuous dependence on the initial data is all that can be expected. Although we use a similar method in proving local well-posedness for \eqref{eq:nonlinear}, the semilinear structure allows us to obtain Lipschitz dependence.

The outline of the proof of Theorem \ref{mainthrm} is as follows. Using a similar argument to Bejenaru-Tataru \cite{BT}, we split the initial data into a low frequency component and a high frequency component. As the low frequency component of the data is essentially stationary on a small time interval, we freeze it at \(t=0\) and rewrite \eqref{eq:nonlinear} as an equation for the evolution of the high frequency component of the form
\begin{equation}\label{intro:hf}(\partial_t+\partial_x^3+a(x)\partial_x^2)v=\tilde F(x,v,v_x,v_{xx})\end{equation}
As the spaces we use are adapted to the unit time interval, we rescale the initial data so that the high frequency component is sufficiently small to solve \eqref{intro:hf} using a perturbative argument on the unit time interval. The Mizohata-type condition \eqref{eq:mizohata} suggests the term \(a(x)\partial_x^2v\) will not be perturbative, so we include this in the principal part. In order to establish estimates for the linear equation
\begin{equation}\label{intro:modairy}(\partial_t+\partial_x^3+a(x)\partial_x^2)v=f\end{equation}
we conjugate the operator by \(e^{-\frac{1}{3}\int_0^xa(y)\,dy}\) and then find approximate solutions to \eqref{intro:modairy} by solving a suitable Airy equation. To complete the proof of Theorem \ref{mainthrm} we use a contraction mapping argument to solve for the time evolution of the high frequency component of the data.

The structure of remainder of the paper is as follows: in Section 2 we define the function spaces used and prove a number of bilinear estimates. In Section 3 we discuss the rescaling and properties of the rescaled initial data. In Section 4 we prove an estimate for the solution to the linear Airy-type equation \eqref{intro:modairy} and in Section 5 we complete the proof of Theorem \ref{mainthrm}.

\vspace{15pt}\begin{remark}
While our result covers the case of the KdV, mKdV and gKdV, it is far from the best known results for these equations and we refer the reader to \cite{LinPon} for a summary of results and references.

However, even in the case of quadratic nonlinearities involving \(u_{xx}\) with which we are primarily concerned, the argument used to prove Theorem \ref{mainthrm} allows us to relax the assumption \(s>\tfrac{9}{2}\) for particular nonlinearities and initial data.

For data with sufficiently small \(\leb H^s\) norm we can prove local well-posedness without having to rescale the initial data. In this case we can use a contraction mapping argument with the linear estimate of Proposition \ref{propn:MMT} and the bilinear and algebra estimates of Proposition \ref{propn:bilests}. The only restrictions on regularity in this case are from the bilinear and algebra estimates and hence we have well-posedness on the unit time interval provided \(s>\sigma_0\) where \(\sigma_0\) is defined as follows.
\vspace{0.3cm}\renewcommand{\arraystretch}{1.5}\begin{center}
\begin{tabular}{ c | c c }
\(\mathbf{\sigma_0}\)&\multicolumn{2}{c}{\(\mathbf{F}\)\textbf{ contains terms of the form}}\\\hline
\(\tfrac{1}{2}\)&\(u^2\)&\\\hline
\(1\)&\(u^{\alpha_0}\)&\(\alpha_0\geq3\)\\\hline
\(\tfrac{3}{2}\)&\(u^{\alpha_0}u_x\)&\\\hline
\(2\)&\(u^{\alpha_0}u_x^{\alpha_1}\)&\(\alpha_0\geq1\)\\\hline
\multirow{2}{*}{\(\tfrac{5}{2}\)}&\(u^{\alpha_0}u_x^{\alpha_1}u_{xx}\)&\(\alpha_0\geq1\)\\
&\(u_x^{\alpha_1}\)&\\\hline
\(3\)&\(u^{\alpha_0}u_x^{\alpha_1}u_{xx}^{\alpha_2}\)&\(\alpha_0\geq1\)\\\hline
\(\tfrac{7}{2}\)&\(u_x^{\alpha_1}u_{xx}^{\alpha_2}\)&\(\alpha_1\geq1\)\\\hline
\(\tfrac{9}{2}\)&\(u_{xx}^{\alpha_2}\)&\\
\end{tabular}\end{center}\vspace{0.3cm}

In the large data case, in order to ensure the high frequency component of the rescaled initial data is small, the scaling of the \(l^1H^s\) spaces (see Proposition \ref{freqbds}) also require that \(s>\lambda+2\) where
\begin{equation}\label{defn:lambda}
\lambda=\max\left\{\frac{\beta_1+2\beta_2-3}{|\beta|-1}:2\leq|\beta|,\,\beta\leq\alpha,\,c_\alpha\neq0\right\}
\end{equation}
and
\[
F(u,u_x,u_{xx})=\sum\limits_{2\leq|\alpha|\leq m} c_\alpha u^{\alpha_0}u_x^{\alpha_1}u_{xx}^{\alpha_2}
\]
In order to get estimates for the nonlinearity (see Proposition \ref{RBDS}) we need to take \(s>\sigma_0^*\) where
\[
\sigma_0^*=\max\left\{\sigma_0(\beta):2\leq|\beta|,\,\beta\leq\alpha,\,c_\alpha\neq0\right\}
\]
and \(\sigma_0(\beta)\) is determined by the table above. So for large data, Theorem \ref{mainthrm} is true for \(s>s_0\) where 
\begin{equation}\label{defn:s0}
s_0=\max\{\sigma_0^*,\lambda+2\}
\end{equation}
\end{remark}\vspace{15pt}


\section{Spaces}\label{sect:ests}

\subsection{Definitions}

For a Sobolev-type space \(U\) we define the \(l^1_{2j}U\) norm by
\[
\|u\|_{l^1_{2j}U}=\sum\limits_{Q\in\mathcal{Q}_{2j}}\|\chi_Qu\|_{U}
\]
and the \(l^\infty_{2j}U\) norm
\[
\|u\|_{l^\infty_{2j}U}=\sup\limits_{Q\in\mathcal{Q}_{2j}}\|\chi_Qu\|_{U}
\]

We define the local energy space \(X\) (see Remark 3.7 of \cite{KPV3}) with norm
\[
\|u\|_X=\sup\limits_{l\geq0}\sup\limits_{Q\in\mathcal{Q}_l}2^{-l/2}\|u\|_{L^2_{t,x}([0,1]\times Q)}
\]
and have the following local smoothing effect for the Airy linear propagator.

\vspace{15pt}\begin{lemma}
If \(f\in L^2(\mathbb{R})\) then
\[
\|e^{-t\partial_x^3}S_jf\|_X\lesssim 2^{-j}\|S_jf\|_{L^2}
\]
\end{lemma}\vspace{15pt}

We look for solutions to the linear equation in the space \(\leb X^s\subset C([0,1],\leb H^s)\) where
\[
\|u\|_{\leb X^s}^2=\sum\limits_{j\geq0}2^{2js}\|S_ju\|^2_{l_{2j}^1X_j}
\]
and
\[
\|u\|_{X_j}=2^j\|u\|_X+\|u\|_{L^\infty_tL^2_x}
\]
We define the atomic space \(Y\) with atoms \(a\) such that there exist \(j\geq 0\) and \(Q\in\mathcal{Q}_j\) with \(\mathrm{supp}\,a\subset[0,1]\times Q\) and \(\|a\|_{L^2_{t,x}([0,1]\times Q)}\lesssim 2^{-j/2}\) with norm given by
\[
\|f\|_Y=\inf\left\{\sum|c_k|:f=\sum c_ka_k,\;a_k\textrm{ atoms}\right\}
\]
We have the duality relation \(Y^*=X\) with respect to the standard \(L^2\) duality (see \cite{MMT} Proposition 2.1). For the inhomogeneous term in the linear equation we use the space \(\leb Y^s\) where
\[
Y_j=2^jY+L_t^1L_x^2
\]
with norm
\[
\|f\|_{Y_j}=\inf\limits_{f=f_1+f_2}\left\{2^{-j}\|f_1\|_Y+\|f_2\|_{L_t^1L_x^2}\right\}
\]
and as above
\[
\|f\|_{\leb Y^s}^2=\sum\limits_{j\geq0}2^{2js}\|S_jf\|^2_{l^1_{2j}Y_j}
\]
We use the Zygmund space \(C^\gamma_*\) with norm
\[
\|u\|_{C^\gamma_*}=\sup_{j\geq0}\,(2^{\gamma j}\|S_ju\|_{L^\infty})\qquad\gamma>0
\]
We have an algebra estimate
\begin{equation}\label{est:Calg}
\|uv\|_{C^\gamma_*}\lesssim\|u\|_{C^\gamma_*}\|v\|_{C^\gamma_*}
\end{equation}
The H\"older space \(C^\gamma\subset C^\gamma_*\) with the estimate
\begin{equation}\label{est:holderzygmund}
\|u\|_{C^\gamma_*}\lesssim \|u\|_{C^\gamma}
\end{equation}
and when \(\gamma\not\in\mathbb{Z}_+\) we have \(C^\gamma=C^\gamma_*\) (see \cite{T} for details).

\subsection{Estimates}

For several estimates we will replace the partition of unity \(\{\chi_Q\}_{Q\in\mathcal{Q}_{2j}}\) by frequency localized versions \(\tilde\chi_Q\), for example taking \linebreak\(\tilde\chi_Q=S_0\chi_Q\), such that each \(\tilde\chi_Q\sim1\) on \(Q\), is rapidly decreasing off \(Q\) and for a Sobolev-type space \(U\)
\[
\sum\limits_{Q\in\mathcal{Q}_{2j}}\|\chi_Qu\|_{U}\sim\sum\limits_{Q\in\mathcal{Q}_{2j}}\|\tilde\chi_Qu\|_{U}
\]
and as a consequence of being frequency localized
\[
\sum\limits_{Q\in\mathcal{Q}_{2j}}\|S_j(\tilde\chi_Qu)\|_U\sim\sum\limits_{Q\in\mathcal{Q}_{2j}}\|\tilde\chi_QS_ju\|_U
\]

Replacing the partition of unity by frequency localized versions we have the following Bernstein-type inequality for \(1\leq p\leq q\leq\infty\)
\begin{equation}\label{est:bernstein}
\|S_ju\|_{l^1_{2j}L^q_x}\lesssim2^{j\left(\tfrac{1}{p}-\tfrac{1}{q}\right)}\|S_ju\|_{l^1_{2j}L^p_x}
\end{equation}

We also note that for any Sobolev-type space \(U\) we can change interval size
\[
\|u\|_{l^1_{2j}U}\lesssim\left\{\begin{array}{ll}2^{2k-2j}\|u\|_{l^1_{2k}U}&\qquad\textrm{for }j\leq k\vspace{3pt}\\\|u\|_{l^1_{2k}U}&\qquad\textrm{for }j>k\end{array}\right.
\]
To see this, if \(j\leq k\) then for \(Q\in\mathcal{Q}_{2j}\) there exists some \(\tilde Q\in\mathcal{Q}_{2k}\) such that \(Q\subset\tilde Q\) and
\[
\|\chi_Qu\|_{U}\lesssim\|\chi_{\tilde Q}u\|_{U}
\]
and each \(\tilde Q\in\mathcal{Q}_{2k}\) is counted \(2^{2k-2j}\) times in this way. If \(j>k\) then for each \(Q\in\mathcal{Q}_{2j}\),
\[
\|\chi_Qu\|_U\lesssim\sum\limits_{\substack{\tilde Q\in\mathcal{Q}_{2k}\\\tilde Q\subset Q}}\|\chi_{\tilde Q}u\|_U
\]
In the case \(U=L^2\) we can improve this to
\begin{equation}\label{est:changer}
\|u\|_{l^1_{2j}L^2}\lesssim2^{k-j}\|u\|_{l^1_{2k}L^2}\qquad\textrm{for }j\leq k
\end{equation}
by writing
\[
\sum\limits_{Q\in\mathcal{Q}_{2j}}\|\chi_Qu\|_{l^1_{2k}L^2}=\sum\limits_{\tilde Q\in\mathcal{Q}_{2k}}\sum\limits_{\substack{Q\in\mathcal{Q}_{2j}\\Q\subset\tilde Q}}\|\chi_Qu\|_{L^2}
\]
and applying the Cauchy-Schwarz inequality to the sum in \(Q\subset\tilde Q\).

We have the following collection of bilinear estimates.

\vspace{15pt}\begin{proposition}\label{propn:bilests} Suppose \(u,v\colon\mathbb{R}\times\mathbb{R}\rightarrow\mathbb{C}\) and \(a,b\colon\mathbb{R}\rightarrow\mathbb{C}\) then

a) Algebra estimates.
\begin{equation}\label{est:Xalg}\|uv\|_{\leb X^s}\lesssim\|u\|_{\leb X^s}\|v\|_{\leb X^s}\qquad s>1\end{equation}
\begin{equation}\label{est:Halg}\|ab\|_{\leb H^s}\lesssim\|a\|_{\leb H^s}\|b\|_{\leb H^s}\qquad s>\tfrac{1}{2}\end{equation}

b) Bilinear estimates I. For \(s>\tfrac{1}{2}\),\begin{equation}\label{est:XXY}\|uv\|_{\leb Y^s}\lesssim\|u\|_{\leb X^\alpha}\|v\|_{\leb X^\beta}\qquad\alpha,\beta\geq s-2\textrm{ and }\alpha+\beta>s+\tfrac{1}{2}\end{equation}\begin{equation}\label{est:HXY}\|au\|_{\leb Y^s}\lesssim \|a\|_{\leb H^\alpha}\|u\|_{\leb X^\beta}\qquad\alpha,\beta\geq s-1\textrm{ and }\alpha+\beta>s+\tfrac{1}{2}\end{equation}

c) Bilinear estimates II. \begin{equation}\label{est:CXX}\|au\|_{\leb X^s}\lesssim \|a\|_{C^\gamma_*}\|u\|_{\leb X^s}\qquad s>\tfrac{1}{2}\textrm{ and }\gamma>s+1\end{equation}\begin{equation}\label{est:CXY}\|au\|_{\leb Y^s}\lesssim \|a\|_{C^\gamma_*}\|u\|_{\leb Y^s}\qquad s>\tfrac{3}{2}\textrm{ and }\gamma>s\end{equation}\begin{equation}\label{est:CHH}\|ab\|_{\leb H^s}\lesssim \|a\|_{C^\gamma_*}\|b\|_{\leb H^s}\qquad s>1\textrm{ and }\gamma>s\end{equation}
\end{proposition}

\begin{proof}Following \cite{MMT} Proposition 3.1 we consider terms of the form \(S_k(S_iuS_jv)\) and the usual Littlewood-Paley trichotomy.

a) For \eqref{est:Xalg} we consider,

\textbf{High-low interactions.} \(|i-k|<4\) and \(j<i-4\). Using the Bernstein inequality \eqref{est:bernstein} we have,
\begin{align*}\|S_k(S_iuS_jv)\|_{l^1_{2k}X_k}&\lesssim\|S_iuS_jv\|_{l^1_{2i}X_i}\\&\lesssim\|S_iu\|_{l^1_{2i}X_i}\|S_jv\|_{L_{t,x}^\infty}\\&\lesssim2^{\tfrac{1}{2}j}\|S_iu\|_{l^1_{2i}X_i}\|S_jv\|_{L^\infty_tL^2_x}\end{align*}
The symmetric low-high interaction is similar.

\textbf{High-high interactions.} \(|i-j|\leq4\) and \(i,j\geq k-4\).
For \(i>k\) we use the Bernstein inequality \eqref{est:bernstein} at frequency \(\sim 2^k\), Cauchy-Schwarz and then change interval size to get
\begin{align*}\|S_k(S_iuS_jv)\|_{l^1_{2k}X_k}&\lesssim 2^{\tfrac{1}{2}k}\|S_iu\|_{l^1_{2k}X_k}\|S_jv\|_{L^\infty_tL^2_x}\\&\lesssim2^{2i-\tfrac{3}{2}k}\|S_iu\|_{l^1_{2i}X_i}\|S_jv\|_{l^1_{2j}X_j}\end{align*}
The result \eqref{est:Xalg} follows from summation.

The argument for \eqref{est:Halg} is identical to \eqref{est:Xalg} except for the improved interval change result \eqref{est:changer} in the high-high interactions that only requires \(s>\tfrac{1}{2}\)

b) We note that for \(l\leq k\)\begin{equation}\label{est:key2Y}\|f\|_{l^1_{2k}Y}\lesssim 2^l\|f\|_{l^1_{2l}L^2_{t,x}}\end{equation}For \eqref{est:XXY} we consider,

\textbf{High-low interactions.} \(|i-k|<4\) and \(j<i-4\).
Using \eqref{est:key2Y} with \(l=j\) followed by Bernstein's inequality \eqref{est:bernstein},
\begin{align*}\|S_iuS_jv\|_{l_{2k}^1Y_k}&\lesssim 2^{-k}\|S_iuS_jv\|_{l^1_{2k}Y}\\
&\lesssim2^{j-k}\|S_iuS_jv\|_{l^1_{2j}L^2_{t,x}}\\
&\lesssim2^{j-k}\|S_iu\|_{l^\infty_{2j}L^2_{t,x}}\|S_jv\|_{l^1_{2j}L^\infty_{t,x}}\\
&\lesssim2^{\tfrac{5}{2}j-2k}\|S_iu\|_{X_i}\|S_jv\|_{l^1_{2j}L_t^\infty L_x^2}
\end{align*}
The symmetric low-high interaction is similar.

\textbf{High-high interactions.} \(|i-j|\leq4\) and \(i,j\geq k-4\).
For \(j>k\) we use \eqref{est:key2Y} with \(l=k\), change interval size and then use Bernstein's inequality \eqref{est:bernstein}
\begin{align*}\|S_k(S_iuS_jv)\|_{l^1_{2k}Y_k}&\lesssim2^{j-k}\|S_k(S_iuS_jv)\|_{l^1_{2j}L^2_{t,x}}\\
&\lesssim2^{j-\tfrac{1}{2}k}\|S_k(S_iuS_jv)\|_{l^1_{2j}L^2_tL^1_x}\\
&\lesssim2^{j-\tfrac{1}{2}k}\|S_iu\|_{L^\infty_tL^2_x}\|S_jv\|_{l_{2j}^1L^2_{t,x}}\\
&\lesssim2^{j-\tfrac{1}{2}k}\|S_iu\|_{X_i}\|S_jv\|_{l^1_{2j}X_j}
\end{align*}
The result follows from summation.

For \eqref{est:HXY} we have,

\textbf{High-low interactions.} \(|i-k|<4\) and \(j<i-4\).
As for \eqref{est:XXY},
\begin{align*}\|S_iaS_ju\|_{l_{2k}^1Y_k}&\lesssim2^{j-k}\|S_iaS_ju\|_{l^1_{2j}L^2_{t,x}}\\
&\lesssim2^{j-k}\|S_ia\|_{L^2_x}\|S_ju\|_{l^1_{2j}L^\infty_{t,x}}\\
&\lesssim2^{\tfrac{3}{2}j-k}\|S_ia\|_{l^1_{2i}L_x^2}\|S_ju\|_{l^1_{2j}L_t^\infty L_x^2}
\end{align*}

\textbf{Low-high interactions.} \(|j-k|<4\) and \(i<j-4\).
For this case we use the \(L^1_tL^2_x\) norm and switch interval size to get
\begin{align*}\|S_iaS_ju\|_{l_{2k}^1Y_k}&\lesssim\|S_iaS_ju\|_{l^1_{2i}L^1_tL^2_x}\\
&\lesssim \|S_ia\|_{l^1_{2i}L_x^\infty}\|S_ju\|_{l^\infty_{2i}L_{t,x}^2}\\
&\lesssim 2^{\tfrac{3}{2}i-j}\|S_ia\|_{l^1_{2i}L_x^2}\|S_ju\|_{X_j}
\end{align*}

\textbf{High-high interactions.} \(|i-j|\leq4\) and \(i,j\geq k-4\).
Identically to \eqref{est:XXY},
\[\|S_k(S_iaS_ju)\|_{l^1_{2k}Y_k}\lesssim2^{j-\tfrac{1}{2}k}\|S_ia\|_{l^1_{2i}L^2_x}\|S_ju\|_{l^1_{2j}X_j}\]

c) For \eqref{est:CXX} we have,

\textbf{High-low interactions.} \(|i-k|<4\) and \(j<i-4\).
Switching interval size we get
\[\|S_k(S_iaS_ju)\|_{l^1_{2k}X_k}\lesssim2^{i-j}\|S_ia\|_{L^\infty}\|S_ju\|_{l^1_{2j}X_j}\]The condition \(\gamma>s+1\) guarantees that the sum in \(k\sim i\) converges.

\textbf{Low-high interactions.} \(|j-k|<4\) and \(i<j-4\).
This case is straightforwards as
\[\|S_k(S_iaS_ju)\|_{l^1_{2k}X_k}\lesssim\|S_ia\|_{L^\infty}\|S_ju\|_{l^1_{2j}X_j}\]

\textbf{High-high interactions.} \(|i-j|\leq 4\) and \(i,j\geq k-4\).
Switching interval size, for \(j>k\) we have
\[\|S_k(S_iaS_ju)\|_{l^1_{2k}X_k}\lesssim2^{2j-2k}\|S_ia\|_{L^\infty}\|S_ju\|_{l^1_{2j}X_j}\]

For the estimate \eqref{est:CXY} the low-high interaction is identical to \eqref{est:CXX}.

\textbf{High-low interactions.} \(|i-k|<4\) and \(j<i-4\).
\[\|S_k(S_iaS_ju)\|_{l^1_{2k}Y_k}\lesssim\|S_ia\|_{L^\infty}\|S_ju\|_{l^1_{2j}Y_j}\]this requires \(\gamma>s\).

\textbf{High-high interactions.} \(|i-j|\leq 4\) and \(i,j\geq k-4\).
\[\|S_k(S_iaS_ju)\|_{l^1_{2k}Y_k}\lesssim2^{3j-3k}\|S_ia\|_{L^\infty}\|S_ju\|_{l^1_{2j}Y_j}\]

The estimate \eqref{est:CHH} is identical to \eqref{est:CXX} for the low-high and high-high interactions and \eqref{est:CXY} for the high-low interactions.
\end{proof}\vspace{15pt}

\section{Scaling}\label{sect:scales}

Given a solution to \eqref{eq:nonlinear}, we rescale using the scaling corresponding to the `worst' monomial nonlinearity in \(F\),
\[u^{(k)}(t,x)=2^{\lambda k}u(2^{-3k}t,2^{-k}x)\]
where \(\lambda\) is as defined in \eqref{defn:lambda}. We have the corresponding rescaled initial data
\[u_0^{(k)}(x)=2^{\lambda k}u_0(2^{-k}x)\]
The rescaling has the effect of sending high frequencies to low frequencies in the sense that
\[S_j(u_0^{(k)})=(S_{j+k}u_0)^{(k)}\]
To make use of this we define the low and high frequency components of the rescaled initial data
\[u_0^{(k)l}=S_0u_0^{(k)}=(S_{\leq k}u_0)^{(k)}\qquad\qquad u_0^{(k)h}=u_0^{(k)}-u_0^{(k)l}\]

As the low frequency component of the rescaled initial data is essentially stationary on the unit time interval we freeze it at \(t=0\) and define
\[v=u^{(k)}-u_0^{(k)l}\qquad\qquad v_0=u_0^{(k)h}\]
We can then rewrite \eqref{eq:nonlinear} as an equation for \(v\) with coefficients depending on \(u_0^{(k)l}\) and the scaling factor \(k\)
\begin{equation}\label{eq:rescalednonlinear}\left\{\begin{array}{l}(\partial_t+\partial_x^3+a(x)\partial_x^2)v=G(x,v,v_x,v_{xx})+L(x,v,v_x)+R(x)\vspace{3pt}\\v(0)=v_0\end{array}\right.\end{equation}
where
\begin{equation}\label{eq:a}a(x)=\sum_{r=0}^2c_r2^{(r-\lambda-1)k}\partial_x^ru_0^{(k)l}\end{equation}
\(G\) is a polynomial in \(v,v_x,v_{xx}\) of degree \(m\) with no constant or linear terms
\begin{equation}\label{eq:g}G(x,v,v_x,v_{xx})=\sum\limits_{2\leq|\beta|\leq|\alpha|\leq m}G_{\alpha,\beta}(x;k)v^{\beta_0}v_x^{\beta_1}v_{xx}^{\beta_2}\end{equation}
\[G_{\alpha,\beta}(x)=c_{\alpha}2^{(\lambda-3-\lambda|\alpha|+\alpha_1+2\alpha_2)k}(u_0^{(k)l})^{\alpha_0-\beta_0}(\partial_xu_0^{(k)l})^{\alpha_1-\beta_1}(\partial_x^2u_0^{(k)l})^{\alpha_2-\beta_2}\]
\(L\) is linear in \(v,v_x\)
\begin{equation}\label{eq:L}L(x,v,v_x)=\sum\limits_{\substack{2\leq|\alpha|\leq m\\|\beta|=1}}L_{\alpha,\beta}(x)v^{\beta_0}v_x^{\beta_1}\end{equation}
\[L_{\alpha,\beta}(x)=c_{\alpha}2^{(\lambda-3-\lambda|\alpha|+\alpha_1+2\alpha_2)k}(u_0^{(k)l})^{\alpha_0-\beta_0}(\partial_xu_0^{(k)l})^{\alpha_1-\beta_1}(\partial_x^2u_0^{(k)l})^{\alpha_2}\]
and \(R\) is an inhomogeneous term
\begin{equation}\label{eq:R}R(x)=\sum\limits_{2\leq|\alpha|\leq m}R_\alpha(x)-\partial_x^3u_0^{(k)l}\end{equation}
\[R_\alpha(x)=c_{\alpha}2^{(\lambda-3-\lambda|\alpha|+\alpha_1+2\alpha_2)k}(u_0^{(k)l})^{\alpha_0}(\partial_xu_0^{(k)l})^{\alpha_1}(\partial_x^2u_0^{(k)l})^{\alpha_2}\]

In order to solve \eqref{eq:rescalednonlinear} we need estimates on the size of the coefficients and initial data. We have the following proposition giving us estimates on the low and high frequency components of the rescaled initial data.

\vspace{15pt}\begin{proposition}\label{freqbds}~

a) For \(r\geq0\), if \(s>r+1\)
\begin{equation}\label{est:lowfreq}\|\partial_x^ru_0^{(k)l}\|_{\leb H^\sigma}\lesssim2^{(\lambda+1-r)k}\|u_0\|_{\leb H^s}\qquad\sigma\geq0\end{equation}
and if \(s\in(\lambda+2,r+1]\)
\begin{equation}\label{est:lowfreq2}\|\partial_x^ru_0^{(k)l}\|_{\leb H^\sigma}\lesssim2^{-\tfrac{1}{2}(s-\lambda-2)k}\|u_0\|_{\leb H^s}\qquad\sigma\geq0\end{equation}

b) For \(r\geq0\) and \(s>r+\tfrac{1}{2}\)
\begin{equation}\label{est:lowfreqcstar}\|\partial_x^ru_0^{(k)l}\|_{C^\gamma_*}\lesssim2^{(\lambda-r)k}\|u_0\|_{\leb H^s}\qquad\gamma>0\end{equation}

c) \begin{equation}\label{est:highfreq}\|u_0^{(k)h}\|_{\leb H^s}\lesssim2^{-(s-\lambda-2) k}\|u_0\|_{\leb H^s}\end{equation}\end{proposition}
\begin{proof}~

a) For \eqref{est:lowfreq} we consider,
\begin{align*}\|\partial_x^ru_0^{(k)l}\|_{\leb H^\sigma}&\sim\sum\limits_{Q\in\mathcal{Q}_0}\|\chi_Q\partial_x^rS_0(u_0^{(k)})\|_{L^2}\\
&\lesssim\sum\limits_{Q\in\mathcal{Q}_0}\|\chi_Q\partial_x^r(S_{\leq k}u_0)^{(k)}\|_{L^2}\\
&\lesssim\sum\limits_{j=0}^k\sum\limits_{Q\in\mathcal{Q}_0}2^{rj+(\lambda-r)k}\|\chi_QS_{j}u_0(2^{-k}x)\|_{L^2}\\
&\lesssim\sum\limits_{j=0}^k\sum\limits_{Q\in\mathcal{Q}_{-k}}2^{rj+(\lambda+\tfrac{1}{2}-r)k}\|\chi_{Q}S_{j}u_0\|_{L^2}\\
&\lesssim\sum\limits_{j=0}^k\sum\limits_{Q\in\mathcal{Q}_{2j}}2^{(r+1)j+(\lambda+1-r)k}\|\chi_{Q}S_{j}u_0\|_{L^2}\\
&\lesssim2^{(\lambda+1-r)k}\|u_0\|_{\leb H^s}\end{align*}
for \(s>r+1\).

If \(s\in(\lambda+2,r+1)\) we replace the final line by
\begin{align*}\sum\limits_{j=0}^k2^{(r+1)j+(\lambda+1-r)k}\|S_{j}u_0\|_{l^1_{2j}L^2}
&\lesssim2^{-(s-\lambda-2)k}\sum\limits_{j=0}^k2^{(r+1-s)(j-k)}2^{sj}\|S_{j}u_0\|_{l^1_{2j}L^2}\\
&\lesssim2^{-(s-\lambda-2)k}\|u_0\|_{\leb H^s}\end{align*}
If \(s=r+1\) then
\begin{align*}\sum\limits_{j=0}^k2^{sj+(\lambda+1-r)k}\|S_{j}u_0\|_{l^1_{2j}L^2}
&\lesssim2^{\tfrac{1}{2}(\lambda+1-r)k}\sum\limits_{j=0}^k2^{-\tfrac{1}{2}(r-1-\lambda)j}2^{sj}\|S_{j}u_0\|_{l^1_{2j}L^2}\\
&\lesssim2^{-\tfrac{1}{2}(s-\lambda-2)k}\|u_0\|_{\leb H^s}\end{align*}

b) For \eqref{est:lowfreqcstar} we have
\begin{align*}\|\partial_x^ru_0^{(k)l}\|_{C^\gamma_*}&\sim\|\partial_x^rS_0(u_0^{(k)})\|_{L^\infty}\\
&\lesssim\sum\limits_{j=0}^k2^{(\lambda-r)k}\|\partial_x^rS_ju_0\|_{L^\infty}\\&\lesssim\sum\limits_{j=0}^k2^{(\lambda-r)k+(\tfrac{1}{2}+r)j}\|S_ju_0\|_{L^2}\\&\lesssim2^{(\lambda-r)k}\|u_0\|_{\leb H^s}\end{align*}
for \(s>r+\tfrac{1}{2}\).

c) For \eqref{est:highfreq} we have that
\begin{align*}\|S_j(u_0^{(k)})\|_{l^1_{2j}L^2}&\lesssim\sum\limits_{Q\in\mathcal{Q}_{2j}}2^{\lambda k}\|\chi_QS_{j+k}u_0(2^{-k}x)\|_{L^2}\\
&\lesssim\sum\limits_{Q\in\mathcal{Q}_{2j-k}}2^{(\lambda+\tfrac{1}{2})k}\|\chi_QS_{j+k}u_0\|_{L^2}\\
&\lesssim\sum\limits_{Q\in\mathcal{Q}_{2(j+k)}}2^{(\lambda+2)k}\|\chi_QS_{j+k}u_0\|_{L^2}\end{align*}
and hence we have
\begin{align*}\|u_0^{(k)h}\|^2_{\leb H^s}&\sim\sum\limits_{j=1}^\infty2^{2js}\|S_ju_0^{(k)}\|^2_{l^1_{2j}L^2}\\
&\lesssim\sum\limits_{j=1}^\infty2^{2js}2^{(2\lambda+4)k}\|\chi_QS_{j+k}u_0\|^2_{l^1_{2(j+k)}L^2}\\
&\lesssim2^{-2(s-\lambda-2)k}\sum\limits_{j=1}^\infty2^{2(j+k)s}\|\chi_QS_{j+k}u_0\|^2_{l^1_{2(j+k)}L^2}\end{align*}\end{proof}
\vspace{15pt}

Due to the Mizohata-type condition \eqref{eq:mizohata} we do not expect to be able to treat the term \(a\partial_x^2v\) in \eqref{eq:rescalednonlinear} as a perturbation of the linear Airy operator and hence include it in the principal part of the equation. From Proposition \ref{freqbds} (a) and (b) we then have the following estimates for the coefficient \(a\).

\vspace{15pt}\begin{proposition}\label{aBDS}Let \(s>\lambda+2\) then,
\begin{equation}\label{est:aest}\|a(x)\|_{\leb H^\sigma}\lesssim\|u_0\|_{\leb H^s}\qquad\sigma\geq0\end{equation}
\begin{equation}\label{est:Caest}\|a(x)\|_{C^\gamma_*}\lesssim2^{-k}\|u_0\|_{\leb H^s}\qquad\gamma>0\end{equation}
For \(r=1,2\) there exists some \(\delta=\delta(s,F)>0\) such that, \begin{equation}\label{est:daest}\|\partial_x^ra(x;k)\|_{\leb H^\sigma}\lesssim2^{-\delta k}\|u_0\|_{\leb H^s}\qquad\sigma\geq0\end{equation}
\end{proposition}\vspace{15pt}

\section{Linear Estimates}\label{sect:lin}

 We consider the linear equation
\begin{equation}\label{eq:modairy}\left\{\begin{array}{l}(\partial_t+\partial_x^3+a\partial_x^2)v=f\vspace{3pt}\\v(0)=v_0\end{array}\right.\end{equation}
where \(a\colon\mathbb{R}\rightarrow\mathbb{C}\). We aim to prove the following result.

\vspace{15pt}\begin{proposition}\label{result:est}Let \(s>s_0\) where \(s_0\) is as defined in \eqref{defn:s0} and suppose \(a\) satisfies \eqref{est:aest}--\,\eqref{est:daest}. Then for \(k>0\) sufficiently large the equation \eqref{eq:modairy} is locally well-posed in \(\leb H^s\) on the unit time interval \([0,1]\) and the solution satisfies the estimate\begin{equation}\label{est:modairysoln}\|v\|_{\leb X^s}\lesssim C(\|u_0\|_{\leb H^s})(\|v_0\|_{\leb H^s}+\|f\|_{\leb Y^s})\end{equation}\end{proposition}\vspace{15pt}

To find a solution we conjugate the linear operator by \(e^{-\tfrac{1}{3}\int_0^xa(y)\,dy}\). A calculation gives
\[e^{\tfrac{1}{3}\int_0^xa\,dy}(\partial_t+\partial_x^3+a\partial_x^2)e^{-\tfrac{1}{3}\int_0^xa\,dy}w=(\partial_t+\partial_x^3)w-(a_x+\tfrac{1}{3}a^2)w_x+(\tfrac{2}{27}a^3-\tfrac{1}{3}a_{xx})w\]
So if \(w\) solves \begin{equation}\label{eq:veqn}\left\{\begin{array}{rcl}(\partial_t+\partial_x^3)w=&\!\!\!\!g\!\!\!\!&=e^{\tfrac{1}{3}\int_0^xa\,dy}f\vspace{3pt}\\w(0)=&\!\!\!\!w_0\!\!\!\!&=e^{\tfrac{1}{3}\int_0^xa\,dy}v_0\end{array}\right.\end{equation}we expect an approximate solution to \eqref{eq:modairy} to be given by\[v=e^{-\tfrac{1}{3}\int_0^xa\,dy}w\]The sense in which this is an approximate solution is summarized in the following result.

\vspace{15pt}\begin{lemma}\label{modairy}Let \(s>s_0\) and \(a\) satisfy \eqref{est:aest}--\,\eqref{est:daest}.\\Suppose \(v=e^{-\tfrac{1}{3}\int_0^xa\,dy}w\) where \(w\) solves \eqref{eq:veqn} then \begin{equation}\label{est:approxu}\|v\|_{\leb X^s}\leq C(\|u_0\|_{\leb H^s})(\|v_0\|_{\leb H^s}+\|f\|_{\leb Y^s})\end{equation} and there exists some \(\delta=\delta(s,F)>0\) such that the error satisfies the estimate\begin{equation}\label{est:rhserror}\|f-(\partial_t+\partial_x^3+a\partial_x^2)v\|_{\leb Y^s}\leq 2^{-\delta k}C(\|u_0\|_{\leb H^s})(\|v_0\|_{\leb H^s}+\|f\|_{\leb Y^s})\end{equation}\end{lemma}\vspace{15pt}

In order to prove this we start with the following result based on the argument of \cite{MMT} Proposition 4.1.

\vspace{15pt}\begin{proposition}\label{propn:MMT}If \(w\) solves \eqref{eq:veqn} then for any \(s\geq0\)\begin{equation}\label{est:airyest}\|w\|_{\leb X^s}\lesssim\|w_0\|_{\leb H^s}+\|g\|_{\leb Y^s}\end{equation}\end{proposition}
\begin{proof}Localizing at frequency \(2^j\)
\[\left\{\begin{array}{l}(\partial_t+\partial_x^3)w_j=g_j\\w_j(0)=w_{0j}\end{array}\right.\]

Similarly to \cite{MMT} Proposition 4.2, we have the energy estimate
\begin{equation}\label{est:locnrg}\|w_j\|_{L^\infty L^2}^2\lesssim \|w_{0j}\|^2_{L^2}+\|w_j\|_{X_j}\|g_j\|_{Y_j}\end{equation}
If we can prove a local energy decay estimate of the form
\begin{equation}\label{est:locnrgdecay}2^{2j}\|w_j\|^2_X\lesssim\|w_{j}\|_{L^\infty L^2}^2+\|w_j\|_{X_j}\|g_j\|_{Y_j}\end{equation}
then combining these we have the estimate
\begin{equation}\label{est:nosum}\|w_j\|_{X_j}^2\lesssim\|w_{0j}\|_{L^2}^2+\|g_j\|_{Y_j}^2\end{equation}

In order to prove \eqref{est:locnrgdecay} it is enough to show that for any \(l\leq 2j\) and any \(Q\in\mathcal{Q}_{l}\)
\begin{equation}\label{est:locnrgdecay2}2^{2j-l}\|w_j\|^2_{L^2_{t,x}([0,1]\times Q)}\lesssim\|w_{j}\|_{L^\infty L^2}^2+\|w_j\|_{X_j}\|g_j\|_{Y_j}
\end{equation}
This is due to the fact that whenever \(l>2j\) we can cut each \(Q\in\mathcal{Q}_l\) into \(2^{l-2j}\) intervals in \(\mathcal{Q}_{2j}\) and hence for \(l>2j\)
\begin{align*}\sup\limits_{Q\in\mathcal{Q}_l}2^{2j-l}\|w_j\|^2_{L^2_{t,x}([0,1]\times Q)}&\lesssim\sup\limits_{Q\in\mathcal{Q}_l}\sup\limits_{\substack{\tilde Q\in\mathcal{Q}_{2j}\\
\tilde Q\subset Q}}\|w_j\|^2_{L^2_{t,x}([0,1]\times \tilde Q)}\\
&\lesssim\sup\limits_{0\leq i\leq 2j}\sup\limits_{\tilde Q\in\mathcal{Q}_i}2^{2j-i}\|w_j\|^2_{L^2_{t,x}([0,1]\times \tilde Q)}\end{align*}

To prove \eqref{est:locnrgdecay2}, for each \(l\leq2j\) and \(Q\in\mathcal{Q}_l\) we aim to construct a self-adjoint Fourier multiplier \(\mathcal{M}\) so that
\begin{align*}&(\mathrm{M}1)\quad\|\mathcal{M}w_j\|_{L^2_x}\lesssim\|w_j\|_{L^2_x}\\
&(\mathrm{M}2)\quad\|\mathcal{M}w_j\|_X\lesssim\|w_j\|_X\\
&(\mathrm{M}3)\quad\langle[\partial_x^3,\mathcal{M}]w_j,w_j\rangle_{L^2_{t,x}}\gtrsim 2^{2j-l}\|w_j\|^2_{L^2_{t,x}([0,1]\times Q)}-O(\|w_j\|^2_{L^2_{t,x}})\end{align*}
The estimate \eqref{est:locnrgdecay2} then follows from
\[\frac{d}{dt}\langle w_j,\mathcal{M}w_j\rangle=2\,\re\langle(\partial_t+\partial_x^3)w_j,\mathcal{M}w_j\rangle+\langle[\partial_x^3,\mathcal{M}]w_j,w_j\rangle\]

By translation invariance we assume \(Q=[-2^{l-1},2^{l-1}]\) and define\[\mathcal{M}w=m_lw\]
where \(m_l(x)=m(2^{-l}x)\) and \(m'(x)=-\psi^2(x)\) for some real-valued \(\psi\in\mathcal{S}\) localized at frequency \(\lesssim 1\). We choose \(\psi\) such that \(m\) is bounded and decreasing and \(\psi\sim 1\) on \(|x|\leq \tfrac{1}{2}\). Clearly \((\mathrm{M}1)\) and \((\mathrm{M}2)\) follow from this choice of \(\mathcal{M}\).

Integrating by parts we have
\begin{align*}\langle[\partial_x^3,m_l]w_j,w_j\rangle&=\langle(\partial_x^3m_l)w_j,w_j\rangle-3\langle \partial_xm_l\partial_xw_j,\partial_xw_j\rangle\\
&=\langle(\partial_x^3m_l)w_j,w_j\rangle+3\cdot2^{-l}\langle\psi^2(2^{-l}x)\partial_xw_j,\partial_xw_j\rangle\end{align*}
Using the properties of \(\psi\) and the frequency localization of \(w_j\) we get \eqref{est:locnrgdecay2}.

To prove \eqref{est:airyest} we start by taking \(Q\) at scale \(M2^{2j}\) for some \(M\) and consider
\[(\partial_t+\partial_x^3)(\chi_Qw_j)=g_j\chi_Q+[\partial_x^3,\chi_Q]w_j\]
We can replace the \(\chi_Q\) with frequency localized versions as before and assume that they are smooth on scale \(M2^{2j}\) so \(|\partial_x^n\chi_Q|\lesssim_n(M2^{2j})^{-n}\)

The estimate \eqref{est:locnrg} gives us that
\begin{equation}\label{ugh}\sum\limits_Q\|\chi_Qw_j\|_{X_j}\lesssim\sum\limits_Q\left\{\|\chi_Qw_{0j}\|_{L^2}+\|\chi_Qg_j\|_{Y_j}\right\}+\sum\limits_Q\|[\partial_x^3,\chi_Q]w_j\|_{Y_j}\end{equation} Using the above bounds on the derivatives of \(\chi_Q\) we get
\[\sum\limits_Q\|[\partial_x^3,\chi_Q]w_j\|_{L_t^1L_x^2}\lesssim M^{-1}\sum\limits_Q\|\chi_Qw_j\|_{L^\infty_tL^2_x}\]
so taking \(M\) sufficiently large (and \(j\)-independent) we can absorb the last term in \eqref{ugh} into the left-hand side.

To make the transition from \(M2^{2j}\) scale to \(2^{2j}\) scale we use an identical argument to changing interval size to show that for any \(h=h(x)\)\[\sum\limits_{\tilde Q\in M\mathcal{Q}_{2j}}\|\chi_{\tilde Q} h\|_{L^2_x}\sim\sum\limits_{Q\in\mathcal{Q}_{2j}}\|\chi_Q h\|_{L^2_x}\]
\end{proof}\vspace{15pt}

From the bilinear estimate \eqref{est:CXX}, for \(\gamma>s+1\) we have\[\|v\|_{\leb X^s}\lesssim\|e^{-\tfrac{1}{3}\int_0^xa\,dy}\|_{C^\gamma_*}\|w\|_{\leb X^s}\]So from \eqref{est:airyest}, \eqref{est:CXY} and \eqref{est:CHH} we have \[\|v\|_{\leb X^s}\lesssim\|e^{-\tfrac{1}{3}\int_0^xa\,dy}\|_{C^\gamma_*}\|e^{\tfrac{1}{3}\int_0^xa\,dy}\|_{C^\gamma_*}(\|v_0\|_{\leb H^s}+\|f\|_{\leb Y^s})\]
Taking \(\gamma\in\mathbb{Z}_+\) and using \eqref{est:holderzygmund} we have
\[\|e^{\pm\tfrac{1}{3}\int_0^xa\,dy}\|_{C^\gamma_*}\lesssim\|e^{\pm\tfrac{1}{3}\int_0^xa\,dy}\|_{C^\gamma}\lesssim e^{\tfrac{1}{3}\|a\|_{L^1}}\langle\|a\|_{C^{\gamma-1}}\rangle^\gamma\]
Choosing \(\sigma>\gamma-1/2\), from Sobolev imbedding and \eqref{est:aest} we have\[\|a\|_{C^{\gamma-1}}\lesssim\|a\|_{H^{\sigma}}\lesssim\|u_0\|_{\leb H^s}\]and as \(s>1\) whenever \(a\not\equiv0\)
\[\|a\|_{L^1}\lesssim\|a\|_{\leb H^s}\lesssim\|u_0\|_{\leb H^s}\]
To prove \eqref{est:rhserror} we write\[\|f-(\partial_t+\partial_x^3+a\partial_x^2)v\|_{\leb Y^s}=\|(a_x+\tfrac{1}{3}a^2)v_x+(\tfrac{1}{3}a_{xx}+\tfrac{1}{3}aa_x+\tfrac{1}{27}a^3)v\|_{\leb Y^s}\]
We can then apply the bilinear estimate \eqref{est:HXY} to get
\begin{align*}
\|(a_x&+\tfrac{1}{3}a^2)v_x+(\tfrac{1}{3}a_{xx}+\tfrac{1}{3}aa_x+\tfrac{1}{27}a^3)v\|_{\leb Y^s}\\
&\lesssim(\|a_x\|_{\leb H^s}+\|a^2\|_{\leb H^s}+\|a_{xx}\|_{\leb H^s}+\|aa_x\|_{\leb H^s}+\|a^3\|_{\leb H^s})\|v\|_{\leb X^s}
\end{align*}
Using the algebra estimate \eqref{est:Halg} we have\[\|aa_x\|_{\leb H^s}\lesssim\|a_x\|_{\leb H^s}\|a\|_{\leb H^s}\] Using the bilinear estimate \eqref{est:CHH}, for \(\gamma>s\) we have
\[
\|a^2\|_{\leb H^s}\lesssim\|a\|_{C^\gamma_*}\|a\|_{\leb H^s}
\]
The estimate \eqref{est:rhserror} then follows from \eqref{est:aest}--\eqref{est:daest}. This completes the proof of Lemma \ref{modairy}.\vspace{15pt}

We can now construct a solution to \eqref{eq:modairy} by iteration. Let \linebreak\(v^{(n)}\) be the approximate solution to \[\left\{\begin{array}{l}(\partial_t+\partial_x^3+a\partial_x^2)v^{(n)}=f^{(n)}\vspace{3pt}\\v^{(n)}(0)=v^{(n)}_0\end{array}\right.\] constructed in Lemma \ref{modairy}, where \(f^{(0)}=f\), \(v_0^{(0)}=v_0\) and for \(n\geq0\)\[\begin{array}{rl}f^{(n+1)}&\!\!\!\!=f^{(n)}-(\partial_t+\partial_x^3+a\partial_x^2)v^{(n)}\\v^{(n+1)}&\!\!\!\!=0\end{array}\]
Then for \(k\) sufficiently large \[v=\sum\limits_{n=0}^\infty v^{(n)}\] converges in \(\leb X^s\) to a solution of \eqref{eq:modairy} satisfying the estimate \eqref{est:modairysoln}.

To prove uniqueness we consider the solution to
\begin{equation}\label{eq:modairyaaaa}\left\{\begin{array}{l}(\partial_t+\partial_x^3+a\partial_x^2)v=0\vspace{3pt}\\v(0)=0\end{array}\right.\end{equation}
Taking \(w=e^{\tfrac{1}{3}\int_0^xa(y)\,dy}v\), we have
\[(\partial_t+\partial_x^3)w=(a_x+\tfrac{1}{3}a^2)w_x-(\tfrac{2}{27}a^3-\tfrac{1}{3}a_{xx})w\]
and as in the proof of Lemma \ref{modairy}
\[\|w\|_{l^1X^s}\lesssim2^{-\delta k}C(\|u_0\|_{\leb H^s})\|w\|_{l^1X^s}\]
so for sufficiently large \(k\) we must have \(w=0\).

\section{Proof of Theorem \ref{mainthrm}}\label{sect:final}


We have the following estimates for the terms on the right-hand side of \eqref{eq:rescalednonlinear}.

\vspace{15pt}\begin{proposition}\label{RBDS}For \(s>s_0\) where \(s_0\) is as defined in \eqref{defn:s0},
\begin{equation}\label{est:gest}\|G(x,v,v_x,v_{xx})\|_{\leb Y^s}\lesssim C(\|u_0\|_{\leb H^s})\|v\|^2_{\leb X^s}\langle\|v\|_{\leb X^s}\rangle^{m-2}\end{equation}
There exists \(\delta_1=\delta_1(s,F)\in(0,1]\) such that,
\begin{equation}\label{est:lest}\|L(x,v,v_x)\|_{\leb Y^s}\lesssim2^{-\delta_1 k}C(\|u_0\|_{\leb H^s})\|v\|_{\leb X^s}\end{equation}
\begin{equation}\label{est:rest}\|R(x)\|_{\leb Y^s}\lesssim2^{-\delta_1 k}C(\|u_0\|_{\leb H^s})\|u_0\|_{\leb H^s}\end{equation}\end{proposition}
\begin{proof} Using the bilinear estimates \eqref{est:XXY} and \eqref{est:CHH} and the algebra estimate \eqref{est:Xalg}, for \(\gamma>s+1\) we have\[\|G(v)\|_{\leb Y^s}\lesssim\sum\limits_{2\leq|\beta|\leq|\alpha|\leq m}\|G_{\alpha,\beta}\|_{C^\gamma_*}\|v\|_{\leb X^s}^{|\beta|}\]
Using the algebra estimate \eqref{est:Calg} and the low frequency estimate \eqref{est:lowfreqcstar} we have\[\|G_{\alpha,\beta}\|_{C^\gamma_*}\lesssim c_\alpha2^{(\lambda(1-|\beta|)-3+\beta_1+2\beta_2)k}\|u_0\|_{\leb H^s}^{|\alpha|-|\beta|}\lesssim\|u_0\|_{\leb H^s}^{|\alpha|-|\beta|}\]

For \(L\) we use the bilinear estimate \eqref{est:HXY} to get
\[
\|L(v)\|_{\leb Y^s}\lesssim\sum\limits_{\substack{2\leq|\alpha|\leq m\\|\beta|=1}}\|L_{\alpha,\beta}\|_{\leb H^s}\|v\|_{\leb X^s}
\]
We then use \eqref{est:CHH} to estimate exactly one of the low frequency terms in \(L_{\alpha,\beta}\) in \(\leb H^s\) using \eqref{est:lowfreq} and the rest in \(C^\gamma_*\) for \(\gamma>s\) using \eqref{est:Calg} and \eqref{est:lowfreqcstar}. Whenever we can apply \eqref{est:lowfreq} to the term in \(\leb H^s\) we get
\[
\|L_{\alpha,\beta}\|_{\leb H^s}\lesssim2^{(\beta_1-2)k}\|u_0\|_{\leb H^s}^{|\alpha|-1}\lesssim2^{-k}\|u_0\|_{\leb H^s}^{|\alpha|-1}
\]
If we have to apply \eqref{est:lowfreq2} to the term in \(\leb H^s\) we get
\[
\|L_{\alpha,\beta}\|_{\leb H^s}\lesssim2^{-\tfrac{1}{2}(s-\lambda-2)k}\|u_0\|_{\leb H^s}^{|\alpha|-1}
\]

We note that as \(R\) does not depend on \(t\), \[\|R\|_{\leb Y^s}\leq\|R\|_{\leb H^s}\]
As in the estimate of \(L_{\alpha,\beta}\) we can estimate exactly one low frequency term in \(R_\alpha\) in \(\leb H^s\) and the rest in \(C^\gamma_*\) for \(\gamma>s\) so
\[\|R_\alpha\|_{\leb H^s}\lesssim2^{-\min\{\tfrac{1}{2}(s-\lambda-2),2-\lambda\}k}\|u_0\|_{\leb H^s}^{|\alpha|}\]
We use \eqref{est:lowfreq}, \eqref{est:lowfreq2} to get
\[\|\partial_x^3u_0^{(k)l}\|_{\leb H^s}\lesssim2^{-\min\{\tfrac{1}{2}(s-\lambda-2),2-\lambda\}k}\|u_0\|_{\leb H^s}\]\end{proof}\vspace{15pt}

We state the following result which we will use to estimate the differences of the polynomial terms \(a,G,L,R\).
\vspace{15pt}\begin{lemma}L\label{polygen}et \(p(x_1,\dots,x_r)=x_1^{\alpha_1}\dots x_r^{\alpha_r}\) then\[p(u_1,\dots,u_r)-p(v_1,\dots,v_r)=\sum\limits_{j=1}^r(u_j-v_j)q_j\]where \[q_j=\sum\limits_{k=0}^{\alpha_j-1}u_1^{\alpha_1}\dots u_{j-1}^{\alpha_{j-1}}u_j^kv_j^{\alpha_j-1-k}v_{j+1}^{\alpha_{j+1}}\dots v_r^{\alpha_r}\]Further, if \(\mathcal{Y}\) is a Banach space and \(\mathcal{X}_1,\dots,\mathcal{X}_r\) are Banach algebras and we have the multilinear estimate\[\|u_1\dots u_r\|_{\mathcal{Y}}\lesssim\|u_1\|_{\mathcal{X}_1}\dots\|u_r\|_{\mathcal{X}_r}\]then\begin{equation}\label{est:polygen}\|p(u_1,\dots,u_r)-p(v_1,\dots,v_r)\|_{\mathcal{Y}}\leq\sum\limits_{j=1}^nC_j\|u_j-v_j\|_{\mathcal{X}_j}\end{equation} where \[C_j=\sum\limits_{k=0}^{\alpha_j-1}\|u_1\|_{\mathcal{X}_1}^{\alpha_1}\dots\|u_{j-1}\|_{\mathcal{X}_{j-1}}^{\alpha_{j-1}}\|u_j\|_{\mathcal{X}_j}^k\|v_j\|_{\mathcal{X}_j}^{\alpha_j-1-k}\|v_{j+1}\|_{\mathcal{X}_{j+1}}^{\alpha_{j+1}}\dots\|v_r\|_{\mathcal{X}_r}^{\alpha_r}\]\end{lemma}


\vspace{15pt}\subsection{Existence}

We prove the existence of a solution by a fixed point argument. By Proposition \ref{result:est}, for \(k>0\) sufficiently large we can find a solution \(w=\mathcal{T}(v)\in\leb X^s\) to
\begin{equation}\label{eq:iterationscheme}\left\{\begin{array}{l}(\partial_t+\partial_x^3+a\partial_x^2)w=G(v)+L(v)+R\vspace{3pt}\\w(0)=v_0\end{array}\right.\end{equation}where \(a,G,L,R\) are as in \eqref{eq:rescalednonlinear}. Let \(\delta_1>0\) be as in Proposition \ref{RBDS} and let \(\sigma=\tfrac{1}{2}\delta_1\). Then define \[K=\{v\in\leb X^s:\|v\|_{\leb X^s}\leq 2^{-\sigma k}\|u_0\|_{\leb H^s}\}\]

\vspace{15pt}\begin{proposition}For \(k>0\) sufficiently large the map \(\mathcal{T}\colon K\rightarrow K\) is a contraction.
\end{proposition}
\begin{proof}
Suppose \(v\in K\) then from \eqref{est:highfreq}, \eqref{est:modairysoln} and Proposition \ref{RBDS},\begin{align*}\|\mathcal{T}(v)\|_{\leb X^s}&\lesssim C(\|u_0\|_{\leb H^s})(\|v_0\|_{\leb H^s}+\|G(v)\|_{\leb Y^s}+\|L(v)\|_{\leb Y^s}+\|R\|_{\leb Y^s})
\\&\lesssim C(\|u_0\|_{\leb H^s})(\|v_0\|_{\leb H^s}+\|v\|^2_{\leb X^s}\langle\|v\|_{\leb X^s}\rangle^{m-2}+2^{-\delta_1k}\|v\|_{\leb X^s}\\&\quad+2^{-\delta_1 k}\|u_0\|_{\leb H^s})
\\&\lesssim C(\|u_0\|_{\leb H^s})(2^{-(s-\lambda-2)k}+2^{-\delta_1k})\|u_0\|_{\leb H^s}\end{align*}
From the proof of Proposition \ref{RBDS} we have \(\sigma<\delta_1<\min\{1,s-\lambda-2\}\) so for sufficiently large \(k\), \[\|\mathcal{T}(v)\|_{\leb X^s}\leq2^{-\sigma k}\|u_0\|_{\leb H^s}\]

For \(v\jT\in K\) the difference \(w=\mathcal{T}(v\oT)-\mathcal{T}(v\tT)\) satisfies
\[\left\{\begin{array}{l}(\partial_t+\partial_x^3+a\partial_x^2)w=G(v\oT)-G(v\tT)+L(v\oT-v\tT)\vspace{3pt}\\w(0)=0\end{array}\right.\]
From \eqref{est:modairysoln} we have\[\|w\|_{\leb X^s}\lesssim C(\|u_0\|_{\leb H^s})(\|G(v\oT)-G(v\tT)\|_{\leb Y^s}+\|L(v\oT-v\tT)\|_{\leb Y^s})\]
Using Lemma \ref{polygen} with the estimates in Proposition \ref{RBDS} we have
\begin{align*}\|G(v\oT)-G(v\tT)\|_{\leb Y^s}\lesssim 2^{-\sigma k}C(\|u_0\|_{\leb H^s})\|v\oT-v\tT\|_{\leb X^s}\end{align*}
From \eqref{est:lest}\[\|L(w)\|_{\leb Y^s}\lesssim2^{-\sigma k}C(\|u_0\|_{\leb H^s})\|v\oT-v\tT\|_{\leb X^s}\]so\[\|\mathcal{T}(v\oT)-\mathcal{T}(v\tT)\|_{\leb X^s}\lesssim2^{-\sigma k}C(\|u_0\|_{\leb H^s})\|v\oT-v\tT\|_{\leb X^s}\]so \(\mathcal{T}\) is a contraction for sufficiently large \(k\).\end{proof}\vspace{15pt}

From the contraction mapping theorem we have a fixed point of \(\mathcal{T}\) in \(K\) which is a solution of \eqref{eq:rescalednonlinear} satisfying\begin{equation}\label{est:unibdmain}\|v\|_{\leb X^s}\leq 2^{-\sigma k}\|u_0\|_{\leb H^s}\end{equation}After adding the low frequency component of the rescaled initial data and rescaling we get a solution \(u\) to \eqref{eq:nonlinear} in \(C([0,2^{-3k}],\leb H^s)\).

\subsection{Lipschitz dependence on initial data and uniqueness}

Suppose we have two solutions, \(u\oT,u\tT\) to \eqref{eq:nonlinear} with initial data \(u_0\oT,u_0\tT\) respectively. Rescaling both with the same value of \(k\) and subtracting the low frequency components from each we get \(v\oT,v\tT\) satisfying\[\left\{\begin{array}{l}(\partial_t+\partial_x^3+a\jT\partial_x^2)v\jT=G\jT(v\jT)+L\jT(v\jT)+R\jT\vspace{3pt}\\v\jT=v_0\jT\end{array}\right.\qquad j=1,2\]We then have\begin{align}\label{eq:diffdata}(\partial_t+\partial_x^3+a\oT\partial_x^2)(v\oT-v\tT)=&\;G\oT(v\oT)-G\tT(v\tT)\\&+L\oT(v\oT)-L\tT(v\tT)\notag\\&+R\oT-R\tT+(a\tT-a\oT)\partial_x^2v\tT\notag\end{align}
Writing \[G\jT=G((u_0\jT)^{(k)l},\partial_x(u_0\jT)^{(k)l},\partial_x^2(u_0\jT)^{(k)l},v\jT,v\jT_x,v\jT_{xx})\qquad j=1,2\] for a polynomial \(G\) with coefficients depending only on \(\lambda\) we can apply Lemma \ref{polygen} with the estimates in Proposition \ref{RBDS}. From \eqref{est:unibdmain} we have \[\|v\jT\|_{\leb X^s}\lesssim 2^{-\sigma k}\|u_0\jT\|_{\leb H^s}\] so as \(G\) is at least quadratic in \(v\jT,v\jT_x,v\jT_{xx}\) we have
\vspace{1pt}\begin{align*}\|G&\oT(v\oT)-G\tT(v\tT)\|_{\leb Y^s}\\&\lesssim C(\|u_0\oT\|_{\leb H^s},\|u_0\tT\|_{\leb H^s})\left(\|u_0\oT-u_0\tT\|_{\leb H^s}+2^{-\sigma k}\|v\oT-v\tT\|_{\leb X^s}\right)\end{align*}\vspace{1pt}
Similarly,\begin{align*}\|L&\oT(v\oT)-L\tT(v\tT)\|_{\leb Y^s}\\&\lesssim 2^{-\delta_1k}C(\|u_0\oT\|_{\leb H^s},\|u_0\tT\|_{\leb H^s})\left(\|u_0\oT-u_0\tT\|_{\leb H^s}+\|v\oT-v\tT\|_{\leb H^s}\right)\end{align*}\vspace{1pt}\[\|R\oT-R\tT\|_{\leb Y^s}\lesssim2^{-\delta_1 k}C(\|u_0\oT\|_{\leb H^s},\|u_0\tT\|_{\leb H^s})\|u_0\oT-u_0\tT\|_{\leb H^s}\]\vspace{3pt}
For the final term we use \eqref{est:XXY}, \eqref{est:aest} and the fact that for any \(f=f(x)\), \(\|f\|_{\leb X^s}\lesssim \|f\|_{\leb H^{s+1}}\) to get\begin{align*}\|(a\oT-a\tT)\partial_x^2v\tT\|_{\leb Y^s}&\lesssim \|a\oT-a\tT\|_{\leb H^{s+1}}\|v\tT\|_{\leb X^s}\\&\lesssim\|u_0\oT-u_0\tT\|_{\leb H^s}\|u_0\|_{\leb H^s}\end{align*}
Applying \eqref{est:modairysoln} to \eqref{eq:diffdata} we have\begin{align*}\|v&\oT-v\tT\|_{\leb X^s}\\&\lesssim C(\|u_0\oT\|_{\leb H^s},\|u_0\tT\|_{\leb H^s})\left(\|u_0\oT-u_0\tT\|_{\leb H^s}+2^{-\sigma k}\|v\oT-v\tT\|_{\leb X^s}\right)\end{align*}so for \(k\) sufficiently large\[\|v\oT-v\tT\|_{\leb X^s}\lesssim C(\|u_0\oT\|_{\leb H^s},\|u_0\tT\|_{\leb H^s})\|u_0\oT-u_0\tT\|_{\leb H^s}\]

We complete the proof of Theorem \ref{mainthrm} by noting that the constant terms depending on \(\|u_0\|_{\leb H^s}\) are of the form \[C(\|u_0\|_{\leb H^s})=(1+e^{C_0\|u_0\|_{\leb H^s}})p(\|u_0\|_{\leb H^s})\] for \(C_0>0\) and polynomially bounded \(p\). So for a sufficiently large constant \(C>0\) we have well-posedness when \(k\geq C\|u_0\|_{\leb H^s}\).

\section*{Acknowledgements}The author would like to thank his advisor Daniel Tataru for suggesting the problem and several key ideas in the proof. He would also like to thank the referee for pointing out an important missing reference.

\bibliographystyle{amsplain}

\end{document}